\documentclass[11pt]{amsart}

\usepackage{amssymb,amsthm,amsmath}

\newtheorem{theorem}{Theorem}[section]

\newtheorem{cor}[theorem]{Corollary}

\theoremstyle{definition}

\newtheorem{example}[theorem]{Example}

\theoremstyle{remark}

\numberwithin{equation}{section}

\newcommand{\ol}{\overline}

\newcommand{\D}{\delta_{max}}

\def \cC {{\mathcal C}}

\def \cE {{\mathcal E}}
\def \cF {{\mathcal F}}
\def \cG {{\mathcal G}}

\def \cM {{\mathcal M}}

\def \Z {\mathbb Z}

\begin{document}

\title[A constructive solution to OP with a large cycle]{A constructive solution to the Oberwolfach problem with a large cycle}


\date{}

\dedicatory{}

\author[Tommaso.Traetta]{Tommaso Traetta}
\address{DICATAM, Universit\`{a} degli Studi di Brescia, Via Branze 43, 25123 Brescia, Italy.}
\email{tommaso.traetta@unibs.it}

\begin{abstract}
For every $2$-regular graph $F$ of order $v$, the Oberwolfach problem 
$OP(F)$ asks whether there is a $2$-factorization of $K_v$ ($v$ odd) or $K_v$ minus a $1$-factor ($v$ even) into copies of $F$. Posed by Ringel in 1967 and extensively studied ever since, this problem is still open.

In this paper we construct solutions to $OP(F)$ whenever $F$ contains a cycle of length greater than an explicit lower bound.
Our constructions combine the amalgamation-detachment technique  with methods aimed at building $2$-factorizations with an automorphism group having a nearly-regular action on the vertex-set.
\end{abstract}

\maketitle

\section{Introduction}
Let $K_v^*$ denote the complete graph $K_v$ if $v$ is odd, or $K_v$ minus the edges of a $1$-factor $I$ if $v$ is even.
Given any $2$-regular graph $F$ of order $v$, the \emph{Oberwolfach problem} $OP(F)$ asks for a decomposition (i.e., a partition of the edge set) of $K^*_v$ into copies of $F$. It was originally posed by Ringel in 1967 when the order $v$ is odd, and then extended to even orders by Huang, Kotzig, and Rosa \cite{HuKoRo79} in 1979. 
We notice that the even variant can be seen as the \emph{maximum packing version} of the original problem posed by Ringel. Although it has received much attention over the past 55 years, OP remains open. 
Complete and constructive solutions are known when the cycles of $F$ have the same length (see \cite{ASSW, Hoffman Schellenberg 91}),
when $F$ is bipartite (see \cite{Bryant Danziger 11}),
when $F$ has exactly two cycles (see \cite{Traetta 13}),
or for orders belonging to an infinite set of primes (see \cite{Bryant Schar 09}) or twice a prime (see \cite{AlBrHoMaSc 16}).

We refer to \cite[Section VI.12]{Handbook} for further results, up to 2006, concerning the solvability of infinite instances of OP, which however settle only a small fraction of the general problem. We notice, in particular, that no complete solution is known as soon as $F$ has 3 or more cycles. A recent survey on constructive resolution methods that proved successful for solving a large portion of the Oberwolfach problem can be found in \cite{BuDaTr24}. We notice that, opposite to the even case, which concerns 2-factorizations of $K_{2n}-I$, the variant to OP that deals with 2-factorizations of $K_{2n}$ with additional copies of a given one factor $I$ has been considered only recently in \cite{BBBSV19, LeSa19, GR, VPM20}.

Letting $L=\{^{\mu_1}\ell_1, \ldots,\,^{\mu_{u}}\ell_u\}$ be a multiset of integers $\ell_1, \ldots, \ell_u\geq3$, with multiplicity $\mu_1, \ldots, \mu_u$, respectively, we write $F\simeq [L]$ or 
$F\simeq [^{\mu_1}\ell_1, \ldots,\,^{\mu_{u}}\ell_u]$, whenever $F$ is a $2$-regular graph whose list of cycle-lengths is $L$. In this case, we may use the notation
$OP(^{\mu_1}\ell_1, \ldots, \,^{\mu_{u}}\ell_u)$ in place of $OP(F)$. We point out that OP$(F)$ has a  solution whenever $|V(F)|\leq 60$ \cite{DFWMR10, SDTBD}, except when ${F} \in \{[^23],[^43],[4,5],[^23,5]\}$.  

Constructive solutions to $OP(F)$ are given in \cite{BuDaTr22} whenever 
\begin{equation}\label{single-flip graphs}
  F \simeq [x, 2m_1, \ldots, 2m_t, \;^2\ell_1, \ldots, \;^2\ell_u],
\end{equation}
and $x$ is greater than an explicit lower bound. 
Condition \eqref{single-flip graphs} places a constraint on the \emph{cycle structure} of $F$ which cannot contain an odd number of $\ell$-cycles for every odd $\ell\neq x$. 
It is worth mentioning that \cite{BuDaTr22} provides a similar result for the \emph{minimum covering} and 
the \emph{$2$-fold variants} of OP, but in the second case (the $2$-fold variant) there is no restriction on the cycle structure of $F$.

The aim of this paper is to generalize the main result of \cite{BuDaTr22} by dropping restriction \eqref{single-flip graphs}. More precisely, we prove Theorem \ref{main} which employs the following notation: given a list $L$ of positive integers (not necessarily distinct), we write $L_0$ and $L_1$ to represent the multiset of even and odd elements of $L$, respectively. Note that by $|L|$ we mean the size of $L$ as a multiset and let $\max (\varnothing) = 0$.

\begin{theorem}\label{main} 
$OP(y,\ell_1, \ell_2, \ldots, \ell_u)$ has an explicit solution whenever 
\[
    y \geq 3b + 24 b_0 + 28 b_1 + 119,
\] 
where
$b=\sum_{i=1}^u \ell_i$, $b_0 = 2|L_0|\,(\max(L_0) +3)$, $b_1 = 7^{|L_1|-1}(2\max(L_1)+1)$ and
$L=\{\ell_1, \ell_2, \ldots, \ell_u\}$.
\end{theorem}
In other words, Theorem \ref{main} (proved in Section 3) constructs solutions to $OP(F)$ for every arbitrary $2$-regular graph $F$ with a large cycle of length greater than an explicit lower bound, thus taking us one step closer to a complete constructive solution of the Oberwolfach Problem.

The emphasis placed by Theorem \ref{main} on providing `explicit' solutions and an `explicit' lower bound aims to point out the constructive approach used in this paper, which is antithetical to the purely existential ones, such as those based on probabilistic methods that recently allowed 
in \cite{GJKK21} 
to obtain a non-constructive asymptotic solution to the Oberwolfach problem for orders greater than a lower bound that is however unquantified.

Theorem \ref{main} exploits two results, Corollary \ref{H-F:cor} and Theorem \ref{pyramidal+mathing:thm},  obtained via completely different methods, described in Section \ref{preliminaries}. 
Corollary \ref{H-F:cor}, proven in \cite{HiJo01}, extends $(1,2)$-decompositions (see Section \ref{dectofacto}) of $K_m^*$ to $2$-factorizations of $K^*_{n}$ ($m<n$) by making use of the very powerful \emph{amalgamation-detachment technique} introduced by Hilton \cite{Hi84}. Theorem \ref{pyramidal+mathing:thm} provides solutions to $OP(x,\;^2\ell_1, \ldots, \;^2\ell_u)$ satisfying the \emph{matching property} (see Section \ref{pyramidal}). Here, the method used is based on constructing solutions to OP with a \emph{pyramidal automorphism group}.

\section{Preliminaries}\label{preliminaries}
In this section we introduce the basic notions and the preliminary results we need to prove Theorem \ref{main}.

Throughout the paper, all graphs are simple and finite. Given a subgraph $F$ of $G$ (briefly, $F\subset G$), we denote by $G\setminus F$ the graph $G$ minus the edges of $F$.
We refer to the number $\ell$ of edges of a \emph{path} or a \emph{cycle} as their length, and speak of an $\ell$-path ($\ell\geq 1$) or an $\ell$-cycle ($\ell\geq 3$), respectively.
In particular, we denote by $P=\langle p_1, \ldots, p_{\ell}\rangle$ the $\ell$-path whose edges are $\{p_i, p_{i+1}\}$, for $1\leq i \leq \ell-1$ and write $(p_1, \ldots, p_\ell)$ to denote the $\ell$-cycle obtained from $P$ by adding the edge $\{p_1, p_\ell\}$.
A \emph{linear forest} (resp. \emph{$2$-regular} graph) is the vertex-disjoint union of paths (resp. cycles) with at least one vertex of degree $2$. We are hence preventing a linear forest to be a \emph{matching}, that is, the vertex disjoint union of $1$-paths.

A graph $F$, which is not a matching, whose vertices have degree 1 or 2, will be called a \emph{$(1,2)$-graph}. 
Hence, $F$ can be either a linear forest, a $2$-regular graph, or the vertex-disjoint union of two such graphs. 
The list of the cycle-lengths of $F$ is referred to as the \emph{cycle structure} of $F$  and denoted by $cs(F)$. Therefore, $cs(F)=\{^{\mu_1}\ell_1, \ldots,\;^{\mu_u} \ell_u\}$ means that $F$ is the vertex-disjoint union of a (possibly empty) linear forest and $u$ distinct $2$-regular graphs, each containing exactly ${\mu_i}\geq 1$ cycles of length $\ell_i\geq 3$, for $1\leq i\leq u$.
An arbitrary $2$-regular graph with cycle structure $L=\{\ell_1, \ldots, \ell_u\}$ will be denoted by $[L]$ or $[\ell_1, \ldots, \ell_u]$. If $F$ is isomorphic to $[L]$, we write $F\simeq [L]$.

A \emph{factor} of a simple graph $G$ 
is a subgraph $F$ of $G$ such that $V(F) = V(G)$. When $F$ is a matching, $2$-regular, or a $(1,2)$-graph, we speak of
a $1$-factor, 2-factor, or $(1,2)$-factor of $G$, respectively.
We recall that $K_v^*$ is either $K_v$ when $v$ is odd, or $K_v\setminus I$ when $v$ is even, where $I$ is a $1$-factor of $K_v$.

A \emph{decomposition} of a simple graph $G$ is a set $\cG$ of graphs whose edge-sets partition $E(G)$. 
We speak of a $2$-decomposition  or a $(1,2)$-decomposition if each graph in $\cG$ is $2$-regular or a $(1,2)$-graph, respectively. Furthermore,
if all graphs in $\cG$ are also factors of $G$, then we speak of a \emph{factorization}, $2$-factorization or $(1,2)$-factorization of $G$, respectively.\\

The proof of our main result, Theorem \ref{main}, is based on  Corollary \ref{H-F:cor} and Theorem \ref{pyramidal+mathing:thm}, obtained through completely different methods which we describe in the following.

\subsection{Extending $(1,2)$-decompositions to $2$-factorizations}
\label{dectofacto}

Corollary \ref{H-F:cor} is based on the \emph{amalgamation-detachment technique} introduced by Hilton \cite{Hi84} to extend a path decomposition of $K_m$ to a Hamiltonian cycle decomposition of $K_{n}$ 
($m<n$). This constructive method was then used in \cite{HiJo01} to solve $OP(x, \ell, \ldots, \ell)$ for every sufficiently large $x$.

We start by recalling the crucial result in \cite{HiJo01}, that is, Theorem
\ref{H-J}, and then provide a reduced version of it, Corollary \ref{H-F:cor}, stated by using the terminology of graph decompositions. 
Theorem \ref{H-J} requires the basic notions on edge-colored graphs, which we recall in the following.

An edge-coloring of a simple graph $G$ with $t$ colors is a function $\gamma$ mapping $E(G)$ onto a set 
$C=\{c_1, \ldots, c_t\}$ of colors. In this case, we say that $G$ is $t$-edge-colored and for each $i$ we denote by $G(c_i)$ the subgraph induced by the edges colored $c_i$. It is not difficult to see that 
$\cG = \{G(c_1), \ldots, G(c_t)\}$ is a decomposition of $G$. Conversely, any decomposition of $G$ naturally induces an edge coloring of $G$, by assigning distinct colors to the graphs of the decomposition. Given a graph $G'$ containing $G$, an edge-coloring $\gamma'$ of $G'$ is called an extension of $\gamma$, if $\gamma'$ coincides with $\gamma$ over the edges of $G$.

Finally, we denote by $\delta_{max}(G)$ the maximum degree of the vertices of $G$, and recall that a composition of $n$ is a sequence $(s_1,\ldots, s_t)$ of positive integers that sums to $n$.

\begin{theorem}\cite[Theorem 6]{HiJo01}\label{H-J}
  Let $m$ and $n$ be integers, $1\leq m\leq n$, and let $(s_1,\ldots, s_t)$ be a composition of $n-1$, where
  each $s_i\in\{1,2\}$. Let $K_m$ be edge-colored with $t$ colors $c_1, \ldots, c_t$, and let $f_i$ denote the number of edges colored $c_i$. This coloring can be extended to an edge-coloring of $K_n$ where each $K_n(c_i)$ is an $s_i$-factor, and if $s_i=2$ then $K_n(c_i)$ contains exactly one more cycle than $K_m(c_i)$, if and only if the following conditions hold for $1\leq i\leq t$:
\begin{enumerate}
  \item  $f_i \geq s_i \left(m-\frac{n}{2}\right)$,
  \item  $s_in$ is even,
  \item  $\D(K_m(c_i)) \leq s_i$.
\end{enumerate}
\end{theorem}

The following result represents a reduced version of Theorem \ref{H-J}, restated using the terminology of graph decompositions.

\begin{cor}\label{H-F:cor}
Let $\cF = \{F_1, \ldots, F_{b}\}$ be a $(1,2)$-decomposition of $K^*_{2a+\epsilon}$, where $\epsilon\in\{1,2\}$. 
If the following condition holds
\begin{equation}\label{cond:suff}
b\geq 2a  - \min_{i} \textstyle{\frac{|E(F_i)|-\epsilon}{2}},
\end{equation}
then there exists a $2$-factorization 
$\cF^+ = \{F_1^+, \ldots, F_{b}^+\}$ of $K^*_{2b+\epsilon}$ such that
\[
F_i\subset F_i^+\;\;\;\text{and}\;\;\; |cs(F^+_i)| = |cs(F_i)| + 1.
\]
\end{cor}
\begin{proof} Set $m=2a+\epsilon$, $n=2b+\epsilon$, 
$t=b +\epsilon-1$ and let $K^*_m = K_m\setminus I$ when $m$ is even, that is, $\epsilon =2$.
We make use of the $(1,2)$-decomposition $\cF$ and
edge-color $K_m$ with $t$ colors $c_1, \ldots, c_t$ so that \[K_m(c_i) = 
\begin{cases} 
 F_i & \text{if $1\leq i\leq b$},\\
 I   & \text{if $i= b+1$ and $\epsilon = 2$}.
\end{cases}
\]

Each $F_i$ has maximum degree $2$, and we set $s_i = \D(K_m(c_i)) = 2$ for $1\leq i\leq b$, while $s_{b+1} = \D(K_m(c_{b+1})) = 1$ when $\epsilon=2$. Note that $(s_1, \ldots, s_t)$ is a composition of $n$. Therefore, conditions $2$ and $3$ of Theorem \ref{H-J} are satisfied. Finally, let $f_i = |E(F_i)|$, for $1\leq i \leq b$, while $f_{b+1} = |I| = m/2$ when $m$ is even. 
Considering that $|E(F_1)|\leq m=2a+\epsilon$, by the inequality \eqref{cond:suff} we get $a\leq b$ which implies, when $m$ is even, 
that $f_{b+1} = \frac{m}{2}\geq s_{b+1}(m - \frac{n}{2})$ (condition $1$ of Theorem \ref{H-J} for $i=b+1$).
Furthermore, \eqref{cond:suff} is equivalent to saying that
$n\geq 2m - \min_i f_i$ which is in turn equivalent to condition $1$ of Theorem  \ref{H-J} for $1\leq i\leq b$. Therefore, there is an edge-coloring of $K_n$ with the same $t$ colors $c_1, \ldots, c_t$ such that
\begin{enumerate}
\item 
$K_n(c_i)$ is an $s_i$-factor of $K_n$ containing $K_m(c_i)$, for $1\leq i\leq t$;
\item $K_n(c_i)$ contains exactly one more cycle than $K_n(c_i)$, for $1\leq i\leq b$.
\end{enumerate}
Letting $F_i^+ = K_n(c_i)$ and considering that $F_{i}^+$ is a $1$-factor of $K_n$ if and only if $i=b+1$ and $\epsilon=2$, the set
$\cF^+ = \{F_1^+, \ldots, F_b^+\}$ provides the desired $2$-factorization of $K^*_n$.
\end{proof}

\subsection{Solutions to OP satisfying the matching property}
\label{pyramidal}
As mentioned above, the proof of Theorem \ref{main} is also based on 
Theorem \ref{pyramidal+mathing:thm} which constructs solutions to
suitable instances of OP satisfying the matching property ($\cM$) defined below.  
Theorem \ref{pyramidal+mathing:thm} is partly
proven in \cite{BuTr12, BuDaTr22} where the methods used are aimed at building
\emph{pyramidal $2$-factorizations}. More precisely,
a solution to $OP(F)$ is called pyramidal if it has an automorphism group $\Gamma$ fixing $1$ or $2$ vertices 
(according to the parity of $|V(F)|$) and acting sharply transitively on the remaining. Pyramidal solutions can be equivalently described as follows: if $|V(F)| = 2k + \epsilon$ with $\epsilon\in\{1,2\}$, a solution $\cF$ to $OP(F)$ is called $\epsilon$-pyramidal (or just pyramidal) over an additive group $\Gamma$ (not necessarily abelian) of order $2k$, if we can label the vertices of $F$ over $\Gamma\;\cup\;\{\infty_1, \infty_\epsilon\}$ so that $\cF = \{F+\gamma\mid \gamma\in \Gamma\}$, where $F+\gamma$ (the right translate of $F$ by $\gamma$) is the graph obtained from $F$ by replacing each vertex $x\not\in\{\infty_1, \infty_\epsilon\}$ with $x+\gamma$. The group of right translations induced by $\Gamma$ over $V(F) = \Gamma\;\cup\;\{\infty_1, \infty_\epsilon\}$ represents an automorphism group of $\cF$ fixing $\infty_1$ and $\infty_\epsilon$ and acting sharply transitively on the remaining vertices.
We notice that one usually uses the term $1$-rotational in place of $1$-pyramidal. 
A more general description of pyramidal $2$-factorizations is given in
\cite{BoMaRi, BuTr15}, while some recent results showing the effectiveness of the pyramidal approach can be found in \cite{BoBuGaRiTr21, BuRiTr17}.

Theorem \ref{starter}, proven in \cite{BuRi08} in a more general setting, shows how to construct a $1$-rotational solution to $OP(\ell_0, \ell_1,\ldots, \ell_u)$.
Letting $\ell_0$ denote the length of the cycles through $\infty_1$, this solution can be easily extended (see \cite{BuTr15, HuKoRo79}) to a $2$-pyramidal solution of $OP(\ell_0+1, \ell_1,\ldots, \ell_u)$.
Before stating  Theorem \ref{starter}, 
we recall that given a graph $G$ with $V(G)\subset \Gamma \,\cup\, 
\{\infty_1, \infty_2\}$, where $\Gamma$ is any additive group not necessarily abelian, the \emph{list of differences} of $G$, is the multiset 
$\Delta G = \{x-y\mid (x,y)\in \Gamma\times\Gamma, \{x,y\}\in E(G)\}$ of all differences of adjacent vertices of $G$ distinct from $\infty_1$ and $\infty_2$.

\begin{theorem}[\cite{BuRi08}]\label{starter}
Let $F\simeq[\ell_0, \ell_1, \ldots, \ell_u]$ with $V(F) = \Z_{2a}\,\cup\, \{\infty_1\}$.
If $F+a = F$ and $\Delta F \supset \Z_{2a}\setminus\{0\}$, then the set
$\cF = \{F+i\mid i\in\Z_{2a}\}$ is a $1$-rotational solution of $OP(F)$.
\end{theorem}

Theorem \ref{pyramidal+mathing:thm} makes use of a general doubling construction described in \cite{BuTr12}. This construction, when applied to a graceful labeling (defined below) of
suitable $(1,2)$-graphs, allows us to construct $2$-regular graphs satisfying the assumptions of Theorem \ref{starter}, hence pyramidal solutions to OP. This first part of Theorem \ref{pyramidal+mathing:thm} (under assumption 1) is proven in \cite[Theorem 6.4]{BuTr12}. However, we recall here its proof since we need to further show that the pyramidal solutions obtained satisfy the matching property ($\cM$) which will be crucial to prove the main result.

From now on, an arbitrary $(1,2)$-graph with cycle structure $L=\{\ell_1, \ldots,$ $\ell_u\}$ and exactly one path-component of length $k$ will be denoted by $[k\mid L]$ or $[k\mid\ell_1, \ldots, \ell_u]$. If the graph $T$ is isomorphic to $[k\mid L]$, we write $T\simeq [k\mid L]$. Such a graph has been called a \emph{zillion graph} in
\cite{BuDaTr22}.

We recall that a \emph{graceful labeling} of $[k\mid \ell_1, \ldots, \ell_u]$ is a graph $T\simeq[k\mid \ell_1, \ldots, \ell_u]$ 
with vertices in $\Z$, such that
   $V(T)     = \{0, \ldots, a\}$ and 
   $\Delta T = \{\pm1, \ldots, \pm a\}$,
 where $a=k+\sum_{i=1}^u \ell_i$.
Graceful labelings of $[k\mid L]$ are built in \cite{BuDaTr22}
whenever $k>B(L)$, where the lower bound $B(L)$ depends on the cycle structure $L$ and it is defined as follows:
\[
\begin{aligned}
   &\text{$B(L) = 6 b_0 + 7 b_1 + 29$ where}, \\
   &\text{$b_0 = 2|L_0|(\max(L_0) +3)$, 
          and $b_1 = 7^{|L_1|-1}(2\max(L_1)+1)$}.
\end{aligned}
\]

\begin{theorem}[\cite{BuDaTr22}]\label{graceful:thm}
$[k\mid L]$ has a graceful labeling whenever $k\geq B(L)$.
\end{theorem}

We call a \emph{halving} of a $2$-regular graph $G\simeq[x,\;^2\ell_1, \ldots, \;^2\ell_u]$ any subgraph $h(G)$ such that
\begin{equation}\label{halving}
  cs(h(G)) = cs(G\setminus h(G)) = \{\ell_1, \ldots, \ell_u\}.
\end{equation} 
Clearly, $h(G)$ can always be obtained by choosing $u$ cycles of lengths $\ell_1, \ldots, \ell_u$, respectively, and then adding an edge of the $x$-cycle of $G$. Note that both $h(G)$ and $G\setminus h(G)$ are certainly $(1,2)$-graphs when $u\geq 1$.

We say that a solution $\cG$ to $OP(x,\;^2\ell_1, \ldots, \;^2\ell_u)$ satisfies the matching property ($\cM$) if 
\begin{enumerate}
  \item[$(\cM)$] there is a matching $M$ such that
  \begin{enumerate}
    \item $cs(G\setminus M) = \{\ell_1, \ldots, \ell_u\} = cs(h(G')\; \cup\; M)$, and
    \item $h(G')\; \cup\; M$ is a $(1,2)$-graph,
  \end{enumerate}  
  for some distinct graphs $G, G'\in\cG$ and some halving $h(G')$
  of $G'$.
\end{enumerate}

\begin{example} \label{example:1}
Here, we show a solution to $OP(x,{}^23,{}^24)$, with $x\in\{3,4\}$, that satisfies the matching property. 

Set $V(K_{17})=\Z_{16}\,\cup\,\{\infty_1\}$ and consider the $2$-factor $G$ of $K_{17}$ defined as the vertex-disjoint union of the following cycles
\[(\infty_1, 2,10), (3,6,4), (11,14,12), (0,5,1,7), (8,13,9,15).
\]
One can check that $\mathcal{G}=\{G+i\mid 1\leq i\leq 8\}$ is a pyramidal solution to $OP({}^33,{}^24)$. Note that $G=G+8\in\mathcal{G}$; also, setting $G'=G+1\in\mathcal{G}$, we have that,
\[
  G'= (\infty, 3,11) \; \cup\; (4,7,5) \; \cup\; (12,15,13) \; \cup\; (1,6,2,8) \; \cup\;  (9,14,10,0)
\]
and $h(G')=\langle\infty,3\rangle \,\cup\, (12,15,13) \,\cup\, (9,14,10,0)$ is a halving of $G'$.
Taking the matching $M=\{\{\infty, 2\}, \{1,5\}, \{4,6\}\}$, one can check that
\begin{enumerate}
  \item $cs(G\setminus M) = \{3,4\} = cs(h(G') \; \cup\; M)$, and
  \item $h(G') \; \cup\; M$ is a $(1,2)$-graph.
\end{enumerate}
Therefore, $\mathcal{G}$ satisfies the matching property.

To obtain a solution to $OP({}^23,{}^34)$ satisfying the matching property we proceed in a similar way. Set $V(K_{18})=\Z_{16}\,\cup\,\{\infty_1, \infty_2\}$ and consider the $2$-factor $G$ of $K_{18}$ defined below:
\[
G = (\infty_1, 2, \infty_2, 10) \,\cup\, (3,6,4) \,\cup\, (11,14,12) \,\cup\,    
    (0,5,1,7) \,\cup\, (8,13,9,15).
\]
One can check that $\mathcal{G}=\{G+i\mid 1\leq i\leq 8\}$ is a pyramidal solution to $OP({}^23,{}^34)$. 
As before $G=G+8\in\mathcal{G}$, and
letting $G'=G+1\in\mathcal{G}$, we have that
\[
  G'=(\infty, 3, \infty_2, 11) \; \cup\; (4,7,5) \; \cup\; (12,15,13) \; \cup\; (1,6,2,8) \; \cup\;  (9,14,10,0)
\]
and $h(G')=\langle\infty,3\rangle \,\cup\, (12,15,13) \,\cup\, (9,14,10,0)$ is a halving of $G'$. One can check that $\mathcal{G}$ satisfies the matching property with respect to $M$, $G$ and $h(G')$.
\end{example}

We end this section building solutions to
$OP(x, \;^2\ell_1,  \ldots,\;^2 \ell_u)$ that satisfies the matching property.

\begin{theorem}\label{pyramidal+mathing:thm} 
There exists a pyramidal solution to 
$OP(x, \;^2\ell_1,  \ldots,\;^2 \ell_u)$, with $u\geq 1$, that satisfies the matching property when either
\begin{enumerate}
  \item there exists a graceful labeling of $\left[\lfloor\frac{x-3}{2}\rfloor\mid \ell_1, \ldots, \ell_u\right]$, or
  \item $x\geq 2B(L)+3$ with $L=\{\ell_1, \ldots, \ell_u\}$.
\end{enumerate}
\end{theorem}
\begin{proof} Let $x=2k+\epsilon$, with $\epsilon\in\{1,2\}$, and set
  $a = k + \sum_{i=1}^u\ell_i$. Note that $a\geq2$ since $u\geq 1$, 
  and $\lfloor\frac{x-3}{2}\rfloor = k-1$.
  Since Theorem \ref{graceful:thm} guarantees the existence of a graceful labeling of $[k-1\mid \ell_1, \ldots, \ell_u]$ whenever $k-1\geq B(L)$, that is, $x\geq 2B(L)+3$, it is enough to prove the assertion under assumption 1. 
  Therefore, let
$T = P\,\cup\,R$ be a graceful labeling of $[k-1\mid\ell_1, \ldots, \ell_u]$,
where $P$ is a $(k-1)$-path disjoint from the $2$-regular graph $R\simeq[\ell_1, \ldots, \ell_u]$. 
Hence,
 \begin{equation}\label{graceful}
   V(T)     = \{0, \ldots, a-1\}\;\;\;\text{and}\;\;\; 
   \Delta T = \{\pm1, \ldots, \pm (a-1)\}.
 \end{equation}
 Let $p_0$ and $p_1$ denote the end-vertices of $P$, let $\cE$ be the graph 
 with $E(\cE) = \{\{\infty_1, p_0\}, \{\infty_1, p_0+a\}, \{p_1, p_1+a\}\}$, and set
 $\cC = P \cup (P+a)\cup \cE$. Finally, set
 \[
   G = T \,\cup\, (T+a)\, \cup\, \cE = 
       R \,\cup\, (R+a)\, \cup\, \cC
 \]
 By \eqref{graceful}, 
 $V(G) = V(T)\,\cup\,V(T+a)\,\cup\,\{\infty_1\} = \{0,\ldots,2a-1\}\,\cup\,\{\infty_1\}$. Also, recalling that $P$ and $R$ are vertex-disjoint, it follows that $P,R, P+a, R+a$ are vertex-disjoint, as well. Hence, $\cC$ is a $(2k+1)$-cycle, and
 $G\simeq[2k+1, \;^2\ell_1,  \ldots,\;^2 \ell_u]$. 
 
From now on, we consider the vertices of $G$ and its subgraphs modulo $2a$, hence $V(G)=\Z_{2a}\,\cup\,\{\infty_1\}$. By \eqref{graceful} and considering that $a\in\Delta \cE$, we have that 
 $\Delta G \supset \Delta T\,\cup\,\Delta \cE \supset  \Z_{2a}\setminus\{0\}$. Also, by construction, $G + a = G$. Therefore, Theorem \ref{starter} guarantees that
 $\cG = \{G+i\mid i\in\Z_{2a}\}$ is a $1$-rotational solution to 
 $OP(2k+1, \;^2\ell_1,  \ldots,\;^2 \ell_u)$.
 
 We now show that $\cG$ satisfies the matching property. 
 Let $H = Q\,\cup\,R$ be a halving of $G$
 obtained by choosing the $1$-path $Q$ in $\cC$ so that
 \[
   V(Q)=
   \begin{cases}
     \{\infty_1, p_0\} & \text{if $1\leq p_0\leq a-1$},\\
     \{\infty_1, p_0+a\} & \text{if $p_0=0$}.     
   \end{cases}
 \]
 Also, let $M = Q\,\cup\, N$ be a matching of $H$, 
 where $N$ consists of $u$ edges belonging to the $u$ distinct cycles of $R$. 
 Recalling that $V(R)\subset\{0, \ldots, a-1\}$, it is not difficult to see that the matching $N$ of $R$ can be chosen so that $0\not\in V(N)$, that is
 \[
   V(N)\subset\{1, \ldots, a-1\}. 
 \]
 Clearly, $cs(G\setminus M) = \{\ell_1, \ldots, \ell_u\}$. 
 Also, $H' = H+(a+1)$ is a halving of $G'=G+(a+1)=G+1$ which is a graph of $\cG$  distinct from $G$, since $a>1$. Note that $H'$ is the vertex-disjoint union of the 
 $1$-path $Q'=Q+(a+1)$ and the $2$-regular graph $R' = R+(a+1)$, where
 \[
   V(R') \subset \{0\}\,\cup\,\{a+1,\ldots, 2a-1\}.
 \]
 Considering the vertex-sets of $Q, N$ and $R'$, we conclude that $M$ and $R'$ are vertex-disjoint. Also, $M\,\cup\,Q'$ is either a matching or a linear forest. Hence
 $M\,\cup\, H'$ is a $(1,2)$-graph with $cs(M\,\cup\, H') = cs(R') = \{\ell_1, \ldots, \ell_u\}$. This proves that $\cG$ satisfies the matching property with respect to $M$, the two distinct graphs $G$ and $G'$, and the halving $H'$ of $G'$, thus showing the assertion for $\epsilon=1$.
 
 Consider  the matching $I$ containing the edges $e_\infty = \{\infty_1, \infty_2\}$ and $e_i=\{p_1+i, p_1+a+i\}$, for $0\leq i< a$. Clearly,
 $I$ is a $1$-factor of $K_{2a+2}$, with 
$V(K_{2a+2}) = \Z_{2a}\,\cup\, \{\infty_1, \infty_2\}$. 
 Now let $G^*$ be the graph obtained from $G$ by removing the edge  
 $e_0$ and then joining its end-vertices to $\infty_2$, and set
 \[
   \cG^* = \{G^*+i\mid 0\leq i <a\}.
 \]
 Note that $G^*\simeq [2k+2, \;^2\ell_1,  \ldots,\;^2 \ell_u]$.
 Since each $G^*+i$ can be obtained from $G+i$ by (performing the same operation of) inserting $\infty_2$ along the edge $e_i\in E(G+i)$, it follows that
 $\cG^*$ is a $2$-factorization of $K_{2a+2}\setminus I$ into copies of $G^*$.
 Hence, $\cG^*$ is a $2$-pyramidal solution of $OP(2k+2, \;^2\ell_1,  \ldots,\;^2 \ell_u)$. 
 Since the operations performed to obtain $\cG^*$ do not involve edges of $M$, 
 then $\cG^*$ satisfies the matching property, with respect to $M$, the two distinct graphs $G^*$ and $(G')^*$, and the halving $H'$ of $(G')^*$, thus showing the assertion for $\epsilon=2$.
\end{proof}

\section{The proof of Theorem \ref{main}}
The idea behind the proof of the main result (Theorem \ref{main}) turns out to be similar to the one that in \cite{HiJo01} allows us to solve $OP(x, \;^{\mu}\ell)$ whenever $x$ is sufficiently large. We start with a solution of $OP(x,\;^2 \ell_1, \ldots,\;^2 \ell_u)$ of order $m$ (Theorem \ref{pyramidal+mathing:thm}) and we decompose each of its factors into two $(1,2)$-graphs (halvings) whose cycle structure is $\{\ell_1, \ldots, \ell_u\}$. We have thus obtained a decomposition $\cG$ of $K^*_m$ into $(1,2)$-graphs  having the same cycle structure, which by means of Corollary \ref{H-F:cor} extends to a solution of $OP(y, \ell_1, \ldots, \ell_u)$ (for a suitable $y$) whose order is $\equiv 1$ or $2\pmod{4}$ according to the parity of $x$
(Theorem \ref{partialf}.(1)). To deal with the remaining classes of orders we need to suitably break one graph of $\cG$ and redistribute its pieces between the remaining graphs of $\cG$ without altering their cycle structure (Theorem \ref{partialf}.(2)). This can be done whenever the initial solution to $OP(x,\;^2 \ell_1, \ldots,\;^2 \ell_u)$ satisfies the matching property (Theorem \ref{pyramidal+mathing:thm}).

\begin{theorem}\label{partialf} 
Let $\cG$ be a solution to $OP(x, \,^2\ell_1, \ldots, \,^2\ell_u)$ ($u\geq 1$) and let $\epsilon \equiv x \pmod{2}$ with $\epsilon\in\{1,2\}$.
Then $OP(y, \ell_1, \ldots, \ell_u)$ is solvable whenever the following conditions hold:
\begin{enumerate}
  \item $y = 2x + 3\sum_{\beta=1}^u \ell_\beta -\epsilon$, or
  \item  $y = 2x + 3\sum_{\beta=1}^u \ell_\beta -\epsilon - 2$, provided that $\cG$ satisfies $(\cM)$.
\end{enumerate}
\end{theorem}

\begin{proof}  
Set $x = 2k+\epsilon$ ($k\geq1$) and $a = k + \sum_{\beta=1}^u \ell_\beta$ 
($\ell_1, \ldots, \ell_u\geq 3$).
Also, let $\cG=\{G_1, \ldots, G_{a}\}$ be a solution to 
$OP(2k +\epsilon, \,^2\ell_1, \ldots, \,^2\ell_u)$, that is, a $2$-factorization of $K^*_{2a+\epsilon}$ where
$cs(G_\alpha)=\{2k+\epsilon, \,^2\ell_1, \ldots, \,^2\ell_u\}$ for
$1\leq\alpha\leq a$.

For each $\alpha$, let $h(G_\alpha)$ be a halving of $G_\alpha$, and set 
\[
F_{2\alpha-1} = h(G_\alpha)\;\;\;\text{and}\;\;\; 
F_{2\alpha} = G_\alpha\setminus h(G_\alpha).
\]
Clearly, $F_{2\alpha-1}$ and $F_{2\alpha}$ decompose $G_\alpha$, hence $\cF=\{F_{i}\mid 1\leq i \leq 2a\}$ is a $(1,2)$-decomposition of 
$K^*_{2a+\epsilon}$, and by \eqref{halving}, 
each $F_i$ is a $(1,2)$-graph such that
$cs(F_i) =\{\ell_1, \ldots, \ell_u\}$.
Since $|\cF|=2a$, condition \eqref{cond:suff} holds, 
hence Corollary \ref{H-F:cor} guarantees the existence of a $2$-factorization 
$\cF^+=\{F^+_{i}\mid 1\leq i \leq 2a\}$  of $K_{4a+\epsilon}$ where 
\[
\text{$F_i\subset F^+_{i}$, and $|cs(F^+_{i})| = |cs(F_i)|+1$,}
\]
which imply that $cs(F^+_{i}) = cs(F_i) \;\cup\; \{y_i\} = \{\ell_1, \ldots, \ell_u, y\}$, and
\[
y_i=4a+\epsilon-\sum_\beta \ell_\beta = 4k + 3\sum_\beta \ell_\beta + \epsilon
= 2x + 3\sum_\beta \ell_\beta - \epsilon = y,
\]
for $1\leq i\leq 2a$. 
In other words, $\cF^+$ is a solution to
$OP(y, \ell_1, \ldots, \,\ell_u)$, and this proves the first part of the theorem.

Now assume that $\cG$ satisfies the matching property $(\cM)$: without loss of generality, we can assume that there is a matching $M$ of $G_1$ such that
\begin{equation}\label{matching}
\begin{aligned}
\text{$cs(G_1\setminus M) =  \{\ell_1, \ldots, \ell_u\} = cs(F_3\, \cup\, M)$, and}\\
\text{$F_3\, \cup\, M$ is a $(1,2)$-graph.}
\end{aligned}
\end{equation}
Set $\ol{F}_2 = G_1\setminus M$, $\ol{F}_3 = F_3 \, \cup \, M$, 
and $\ol{F}_i = F_i$ for $4\leq i\leq 2a$. 
By \eqref{matching}, and recalling that $F_1$ and $F_2$ decompose $G_1$, and $M\subset G_1$, 
it follows that
\[
\ol{\cF} = \left\{\ol{F}_i\mid 2\leq i\leq 2a\right\} 
         = (\cF \setminus\{F_1, F_2, F_3\}) \; \cup \; \left\{\ol{F}_2, \ol{F}_3\right\}
\] 
is a $(1,2)$-decomposition of $K_{2a+\epsilon}$ into $b=2a-1$ graphs. Considering that each $\ol{F}_i$ contains at least one path-component of length $\geq 1$ and a cycle-component of length $\geq 3$, 
it follows that $|E(\ol{F}_i)|\geq 4$ for $2\leq i\leq 2a$. Therefore,
$\ol{\cF}$ satisfies condition \eqref{cond:suff}, 
hence Corollary \ref{H-F:cor} guarantees the existence of a $2$-factorization 
$\ol{\cF}^+=\left\{\ol{F}^+_{i} \mid 2\leq i \leq 2a\right\}$  
of $K_{2b+\epsilon} = K_{4a+\epsilon-2}$ such that
\[
\ol{F_i}\subset \ol{F}^+_{i},
\;\;\; \text{and}\;\;\; 
|cs(\ol{F}^+_{i})| = |cs(\ol{F}_i)|+1,
\]
for $2\leq i\leq 2a$. Reasoning as before, we conclude that 
$\ol{\cF}^+$ is a solution to
$OP(y, \ell_1, \ldots, \,\ell_u)$ where
\[
   y = 4a+\epsilon-2-\sum_{\beta=1}^u \ell_{\beta} 
     = 4k - 3\sum_{\beta=1}^u \ell_\beta +\epsilon-2
     = 2x - 3\sum_{\beta=1}^u \ell_\beta -\epsilon-2,     
\]
and this completes the proof.
\end{proof}

\begin{example} \label{example:2}
Here, we follow the proof of Theorem \ref{partialf}, and by starting with a solution to $OP(x, {}^2 \ell_1, {}^2 \ell_2)$, with $x\in\{3,4\}$ and $(\ell_1, \ell_2)=(3,4)$, we construct a solution to $OP(y, 3, 4)$, for $y\in\{24, 25,26,27\}$.
        
Let $\epsilon\in\{1,2\}$ and take the solution 
$\mathcal{G} = \{G_\alpha\mid 1\leq \alpha \leq a=8\}$ 
of OP$(2+\epsilon,{}^23, {}^24)$ considered in Example \ref{example:1},
where $G_\alpha=G+\alpha$, and $G$ is the $2$-regular graph defined below, according to the values of $\epsilon$: if $\epsilon=1$, then
\[
G= (\infty_1, 2,10) \,\cup\, (3,6,4) \,\cup\, (11,14,12)
                   \,\cup\, (0,5,1,7) \,\cup\, (8,13,9,15),
\]
otherwise,
\[
 G = (\infty_1, 2, \infty_2, 10) \,\cup\, (3,6,4) \,\cup\, (11,14,12)
                   \,\cup\, (0,5,1,7) \,\cup\, (8,13,9,15).
\]       
Note that $G=G_8\in\mathcal{G}$. For $1\leq \alpha\leq 8$, consider the halving 
\[
h(G_\alpha) = \langle\infty,2+\alpha\rangle \,\cup\, (11+\alpha,14+\alpha,12+\alpha) 
\,\cup\, (8+\alpha,13+\alpha,9+\alpha,15+\alpha)
\] of $G_\alpha$, and set $F_{2\alpha-1} = h(G_\alpha)$ and
$F_{2\alpha} = G_\alpha\setminus h(G_\alpha)$. 
Clearly, 
\[\cF=\{F_{i}\mid 1\leq i \leq 2a=16\}\] 
is a $(1,2)$-decomposition of 
$K^*_{16+\epsilon}$. 
Since $b=|\cF|=16$, condition \eqref{cond:suff} holds, 
hence Corollary \ref{H-F:cor} guarantees the existence of a $2$-factorization 
$\cF^+=\{F^+_{i}\mid 1\leq i \leq 16\}$  of $K_{4a+\epsilon}=K_{32+\epsilon}$ where each
$F_i^+$ is a $2$-regular graph of order $32+\epsilon$ containing $F_i$ and exactly one more cycle than $F_i$, of length say $y_i$. Since 
$cs(F_i) =\{3,4\}$, we have that $y_i=25+\epsilon$ and
$\cF^+$ is a solution to $OP(25+\epsilon, 3,4)$.

It is left to build a solution to $OP(24, 3,4)$ and $OP(25, 3,4)$.
In Example \ref{example:1}, we showed that $\mathcal{G}$ satisfies the matching propery $(\cM)$: more precisely, by taking the matching $M=\{\{\infty, 2\}, \{1,5\}, \{4,6\}\}$ of $G_8=G$ and recalling that each $G_\alpha$ decomposes into $F_{2\alpha-1}$ and $F_{2\alpha}$ (both halvings of $G_\alpha$), one can check that 
\begin{equation}\label{matching}
\begin{aligned}
\text{$cs(G_8\setminus M) =  \{3,4\} = cs(F_1\, \cup\, M)$, and}\\
\text{$F_1\, \cup\, M$ is a $(1,2)$-graph.}
\end{aligned}
\end{equation}
Recall that $a=8$ and let $\ol{\cF}$ be the $(1,2)$-decomposition of $K_{2a+\epsilon}$ obtained from $\mathcal{F}$ by replacing $F_{15}$ and $F_{16}$ with $G_8\setminus M$, and then replacing $F_1$ with $F_1\, \cup\, M$.
By \ref{matching}, we have that $cs(\ol{F})=\{3,4\}$, for every
$\ol{F}\in \ol{\cF}$. Also, note that $b=|\ol{\cF}|=2a-1=15$, and considering that $|E(F_i)|\geq 7$, it follows that
$\ol{\cF}$ satisfies condition \eqref{cond:suff}. Hence,
Corollary \ref{H-F:cor} guarantees the existence of a $2$-factorization 
$\ol{\cF}^+=\left\{\ol{F}^+_{i} \mid 1\leq i \leq 2a-1=15\right\}$  
of $K_{2b+\epsilon} = K_{30+\epsilon}$ where each
$\ol{F_i}^+$ is a $2$-regular graph of order $30+\epsilon$ containing $\ol{F}_i$ and exactly one more cycle than $\ol{F}_i$, of length say $y_i$. Since 
$cs(F_i) =\{3,4\}$, we have that $y_i=23+\epsilon$ and
$\cF^+$ is a solution to $OP(23+\epsilon, 3,4)$.
\end{example}

We are now ready to prove the main result of this paper, which we restate below.\\

\noindent
\textbf{Theorem \ref{main}.} 
\emph{$OP(y,\ell_1, \ell_2, \ldots, \ell_u)$ has an explicit solution whenever 
\[
   y \geq 3b + 24 b_0 + 28 b_1 + 119,
\] 
where
$b=\sum_{i=1}^u \ell_i$, $b_0 = 2|L_0|\,(\max(L_0) +3)$, $b_1 = 7^{|L_1|-1}(2\max(L_1)+1)$ and
$L=\{\ell_1, \ell_2, \ldots, \ell_u\}$.}

\begin{proof} Since $OP(y)$ and $OP(y,\ell)$ are completely solved 
(see, \cite[Section VI.12]{Handbook} and \cite{BuDaTr22}), we can assume $u\geq 2$. Consider the maps $\epsilon: \mathbb{N} \rightarrow \{1,2\}$
and $f: \mathbb{N} \rightarrow \mathbb{N}$ defined as follows: 
    \[\epsilon(y) \equiv y+b \;({\rm mod}\;2),\;\;\;
      f(y)=\frac{y+\epsilon(y) - 3b}{2} + 
      \begin{cases}
         0 &\text{if $\epsilon(y) \equiv y+b \;({\rm mod}\;4$)},\\
         1 &\text{otherwise.}\     
      \end{cases}
    \]
Now set $y_0 = 3b + 24 b_0 + 28 b_1 + 119$.
Considering that $y_0 + b \equiv 3 \pmod{4}$, then
$\epsilon(y_0)=1$, and it is not difficult to check that
\[
\min_{y_0\leq y} f(y) = f(y_0) =  12 b_0 + 14b_1 + 61.
\]
Hence, for every $\ol{y}\geq y_0$, Theorem \ref{pyramidal+mathing:thm} constructs 
a pyramidal solution to $OP(f(\ol{y}),$ $\;^2\ell_1, \ldots, \;^2\ell_u)$ that satisfies the matching property.
Finally, Theorem \ref{partialf} constructs a solution to 
$OP(\ol{y}, \ell_1, \ldots, \ell_u)$.  
\end{proof}

\section{Conclusions}
In this paper we construct solutions to the Oberwolfach problem $OP(F)$
for every $2$-regular graph $F$ with a cycle whose length is greater than an explicit lower bound: Theorem \ref{main}.
This result makes use of two results, Corollary \ref{H-F:cor} and Theorem \ref{pyramidal+mathing:thm}  obtained via completely different methods. Corollary \ref{H-F:cor}, proven in \cite{HiJo01}, extends $(1,2)$-decompositions of $K_m^*$ to $2$-factorizations of $K^*_{n}$ ($m<n$) by making use of the very powerful amalgamation-detachment technique introduced by Hilton \cite{Hi84}. Theorem \ref{pyramidal+mathing:thm} constructs solutions to $OP(x,\;^2\ell_1, \ldots, \;^2\ell_u)$ satisfying the matching property. Here, the method used is based on constructing $2$-factorizations with a pyramidal automorphism group.

The main idea behind the proof of Theorem \ref{main} can be easily generalized as follows. 
We start with a solution $\cF=\{F_1, \ldots, F_a\}$ to $OP(x, \,^{\mu}\ell_1,\ldots, \,^{\mu}\ell_u)$ of order $v$.
Then, we decompose each factor $F_i$ of $\cF$ into $\mu$ $(1,2)$-graphs 
$F_{i,1},\ldots, F_{i,{\mu}}$ with the same cycle structure: $cs(F_{i,j}) = \{\ell_1, \ldots, \ell_u\}$. In other words, we separate out $\mu$ sets of cycles of length $\ell_1, \ldots, \ell_u$ and then add to each set a portion of the $x$-cycle of $F$. We have obtained a $(1,2)$-decomposition $\cG$ of $K_v^*$, and by applying Corollary \ref{H-F:cor}, we  construct a $2$-factorization that solves $OP(y, \,\ell_1,\ldots, \,\ell_u)$ for a specific value of $y$. 
In other words, we have proven the following.

\begin{theorem}
  If $OP(x, \,^{\mu}\ell_1,\ldots, \,^{\mu}\ell_u)$ of order $2w+\epsilon$, with $\epsilon\in\{1,2\}$, has a solution,
then there is a solution to $OP(y, \ell_1,\ldots, \ell_u)$ 
with $y=2w{\mu}+\epsilon - \sum \ell_i$.
\end{theorem}

Note that the order of $OP(y, \ell_1,\ldots, \ell_u)$ is 
$2w{\mu}+\epsilon \equiv \epsilon \pmod{2{\mu}}$. To deal with the remaining classes of orders $({\rm mod}\; 2{\mu})$ it is enough to manipulate the intermediate
$(1,2)$-decomposition $\cG$, by decomposing $i$ graphs in $\cG$ ($1\leq i<{\mu}$) into suitable linear forests that can then be added to the remaining graphs in $\cG$ to form larger $(1,2)$-graphs but with the same initial cycle structure. This way we end up producing a solution to 
\begin{equation}\label{general}
 \text{$OP(y, \ell_1,\ldots, \ell_u)$ 
with $y=2w{\mu}+\epsilon - 2i - \sum \ell_i$.}
\end{equation}
Succeeding to solve \eqref{general},
for every $1\leq i< {\mu}$, would then lead to solve $OP(y,$ $\ell_1,\ldots,\ell_u)$ for every $y>f(x_0, {\mu})$, provided that we can solve $OP(x,$ $\,^{\mu}\ell_1,\ldots, \,^{\mu}\ell_u)$ for every $x>x_0$. Note that $f$ would be an increasing function of both $x_0$ and $\mu$. An explicit value for the lower bound $x_0$ is given in \cite{BuDaTr22} and it grows as ${\mu}$ increases.
Therefore, the best possible lower bound on $y$, based on the results of \cite{BuDaTr22}, can be achieved when ${\mu}=2$. This is the reason why
all constructions in this paper make use of solutions to $OP(x, \,^2\ell_1,\ldots, \,^2\ell_u)$.

We conclude with two tables showing the smallest value $\ol{y}$, given by Theorem \ref{main}, that guarantees the solvability of $OP(y, \ell_1, \ell_2)$ for every $y\geq \ol{y}$, when $3\leq \ell_1  < \ell_2\leq 8$. 
We exclude the cases where $\ell_1=\ell_2$ since for them a better lower bound, that is $\ol{y}=5$, is given in \cite{Traetta 13}. Further partial results can be found in \cite{VaMu16}.\\

\begin{tabular}{c|c|c|c|c|c|c|c|c|c|c|c}
$\ol{y}$&  672& 2299& 774& 3089& 876& 3879&
 790& 1017& 908& 1215& 1026 \\ \hline
$\ell_1$&  3& 3& 3& 3& 3& 3&
 4& 4& 4& 4& 4\\ \hline
$\ell_2$ &  4& 5& 6& 7& 8& 9&
 5& 6& 7& 8& 9 \\ \hline
\end{tabular}

\vspace{5mm}

\begin{tabular}{c|c|c|c|c|c|c|c|c|c|c}
$\ol{y}$ & 
892& 3095& 994& 3885&
1010& 1221 & 1128&
1112& 3891& 
1230 \\ \hline
$\ell_1$ & 
5& 5& 5& 5& 
6& 6& 6& 
7& 7& 
8 \\ \hline
$\ell_2$ & 
6& 7& 8& 9&
7& 8& 9&
8& 9&
9 \\ \hline
\end{tabular}

\vspace{5mm}

An improvement to Theorem \ref{main} containing a lower bound on $y$ that is linear in the remaining $u$ cycle lengths will be given in a paper in preparation \cite{Tr}.

\section*{Acknowledgements}
The author was supported by INdAM-GNSAGA.

\bibliographystyle{amsplain}

\begin{thebibliography}{99}
\bibitem{AlBrHoMaSc 16}
B. Alspach, D. Bryant, D. Horsley, B. Maenhaut, V. Scharaschkin,
On factorisations of complete graphs into circulant graphs and the Oberwolfach problem, 
Ars Math. Contemp. 11(2016), 157--173.

\bibitem{ASSW}
B.\ Alspach, P.J.\ Schellenberg, D.R.\ Stinson, and D.\ Wagner,
The Oberwolfach problem and factors of uniform odd length cycles,
{J. Combin. Theory Ser. A} {52}(1989), 20--43.

\bibitem{BBBSV19}
N. Bolohan, I. Buchanan, A. Burgess, M. \v{S}ajna, R. Van Snick, 
On the spouse-loving variant of the Oberwolfach problem,
{J. Combin. Des.} {27}(2019) 251--260.

\bibitem{BoBuGaRiTr21}
S. Bonvicini, M. Buratti, M. Garonzi, G. Rinaldi, T. Traetta, 
The first families of highly symmetric Kirkman Triple Systems whose orders fill a congruence class, 
Des. Codes Cryptogr. 89(2021), 2725--2757

\bibitem{BoMaRi}
S. Bonvicini, G. Mazzuoccolo, R. Rinaldi,
On 2-factorizations of the complete graph: From the $k$-pyramidal to the universal property, 
J. Combin. Des  17(2009), 211--228.

\bibitem{Bryant Danziger 11}
D. Bryant and P. Danziger.
On bipartite 2-factorizations of $K_n - I$ and the Oberwolfach Problem. 
{J. Graph Theory} {68}(2011), 22--37.

\bibitem{Bryant Schar 09}
D. Bryant and V. Scharaschkin. 
Complete solutions to the Oberwolfach Problem for an infinite set of orders. 
J. Combin. Theory Ser. B {99}(2009), 904--918.

\bibitem{BuRi08}
M. Buratti, G. Rinaldi.
$1$-Rotational $k$-factorizations of the complete graph and new solutions to the Oberwolfach problem. 
J. Combin. Des. 16(2008), 87--100.

\bibitem{BuRiTr17}
M. Buratti, G. Rinaldi, T. Traetta. 
3-pyramidal Steiner triple systems. 
Ars Math. Contemp. 13(2017), 95--106.

\bibitem{BuTr12}
M. Buratti, T. Traetta.
2-starters, graceful labelings and a doubling construction for the Oberwolfach problem. 
J. Combin. Des. 20(2012), 483--503.

\bibitem{BuTr15}
M. Buratti, T. Traetta. 
The structure of 2-pyramidal 2-factorizations. 
Graphs Combin. 31(2015), 523--535.

\bibitem{BuDaTr22}
A.C. Burgess, P. Danziger, T. Traetta. 
On the Oberwolfach problem for single-flip 2-factors via graceful labelings. 
J. Combin. Theory Ser. A 189(2022), 105--611.

\bibitem{BuDaTr24}
A. Burgess, P. Danziger, T. Traetta. 
A survey on constructive methods for the Oberwolfach problem and its variants.
To appear on Fields Institute Communications.

\bibitem{Handbook}
C.J.\ Colbourn and J.H.\ Dinitz, editors.
\newblock {The CRC Handbook of Combinatorial Designs}.
\newblock 2nd ed. CRC Press Series on Discrete Mathematics, Boca Raton, 2007.

\bibitem{DFWMR10},
A. Deza, F. Franek, W. Hua, M. Meszka, A. Rosa.
Solutions to the Oberwolfach Problem for orders 18 to 40.
{J. Math. Combin. Comput.}, {74}(2010), 95--102.

\bibitem{GJKK21}
S. Glock, F. Joos, J. Kim, D. K\"{u}hn, D. Osthus. 
Resolution of the Oberwolfach problem. 
J. Eur. Math. Soc. 23(2021), 2511--2547.

\bibitem{Hi84}
A.J.W. Hilton. 
Hamiltonian decompositions of complete graphs. 
J. Combin. Theory Ser. B 36(1984), 125--134.

\bibitem{HiJo01}
A.J.W. Hilton, M. Johnson. 
Some results on the Oberwolfach problem. 
J. London Math. Soc. (2) 64(2001), 513--522.

\bibitem{Hoffman Schellenberg 91}
D.G. Hoffman and P.J. Schellenberg.
The existence of $C_k$-factorizations of $K_{2n} - F$.
{Discrete Math.} {97}(1991), 243--250.

\bibitem{HuKoRo79}
C. Huang, A. Kotzig, A. Rosa. 
On a variation of the Oberwolfach problem. 
Discrete Math. {27}(1979), 261--277.

\bibitem{LeSa19}
D. Lepine, M. Sajna.
On the honeymoon Oberwolfach problem.
{J. Combin. Des.} {27} (2019), 420--447.

\bibitem{GR}
G. Rinaldi.
The Oberwolfach Problem with loving couples, preprint.

\bibitem{SDTBD}
F. Salassa, G. Dragotto, T. Traetta, M. Buratti, F. Della Croce.
Merging Combinatorial Design and Optimization: the Oberwolfach Problem.
{Australas. J. Combin.} {79}(2021), 141--166.

\bibitem{Traetta 13}
T.\ Traetta.
A complete solution to the two-table Oberwolfach Problems.
J. Combin. Theory Ser. A {120}(2013), 984--997.

\bibitem{Tr}
T.\ Traetta.
A linear lower bound for the solvability of the Oberwolfach problem.
In preparation.

\bibitem{VaMu16}
A. S. Vadivu, A. Muthusamy.
Note on three table Oberwolfach problem.
Electron. Notes Discrete Math. 53(2016), 97--112.

\bibitem{VPM20}
A. Vadivu, L. Panneerselvam, A. Muthusamy.
Solution to the outstanding case of the spouse-loving variant of the Oberwolfach problem with uniform cycle length,
{J. Combin. Des.} {29}(2020), 114--124.
\end{thebibliography}

\end{document}